\documentclass[12pt,a4paper]{article}

\usepackage[latin2]{inputenc}
\usepackage{amssymb}
\usepackage{paralist}
\usepackage{amsmath}
\usepackage{amsfonts}
\usepackage{amsthm}
\usepackage{fullpage}
\usepackage{enumerate}
\newtheorem*{thm1}{Theorem 1}

\newtheorem*{lem1}{Lemma 1}
\newtheorem*{lem2}{Lemma 2}
\newtheorem*{lem3}{Lemma 3}
\newtheorem*{lem4}{Lemma 4}
\newtheorem*{lem5}{Lemma 5}
\newtheorem*{lem6}{Lemma 6}

\begin{document}

\title{On Sidon sets which are asymptotic bases of order 4}
\author{Sándor Z. Kiss \thanks{Institute of Mathematics, Budapest University of Technology and Economics, H-1529 B.O. Box, Hungary; Computer and Automation Research Institute of the Hungarian Academy of Sciences, Budapest H-1111, L\'agym\'anyosi street 11; kisspest@cs.elte.hu; This author was supported by the OTKA Grant No. K77476 and No. NK105645.}, Eszter Rozgonyi \thanks{Institute of Mathematics, Budapest University of Technology and Economics, H-1529 B.O. Box, Hungary, reszti@math.bme.hu, The work reported in the paper has been developed in the framework of the project "Talent care and cultivation in the scientific workshops of BME" project. This project is supported by the grant TÁMOP - 4.2.2.B-10/1--2010-0009.}, Csaba Sándor \thanks{Institute of Mathematics, Budapest University of Technology and Economics, H-1529 B.O. Box, Hungary, csandor@math.bme.hu, This author was supported by the OTKA Grant No. K81658.}
}
\date{}
\maketitle

\begin{abstract}
Let $h \geq 2$ be an integer. We say that a set $\mathcal{A}$ of positive integers is an asymptotic basis of order $h$ if every large enough positive integer can be represented as the sum of $h$ terms from $\mathcal{A}$. A set of positive integers $\mathcal{A}$ is called a Sidon set if all the sums $a+b$ with $a,b \in \mathcal{A}$, $a \leq b$ are distinct. In this paper we prove the existence of Sidon set $\mathcal{A}$ which is an asymptotic basis of order $4$ by using probabilistic methods.
\end{abstract}
\textit{2000 AMS \ Mathematics subject classification number}: primary: 11B13, secondary: 11B75.


\section{Introduction}
Let $\mathbb{N}$ denote the set of nonnegative integers. Let $\mathcal{A} = \{a_{1}, a_{2},\ldots{}\}$ $(a_{1} < a_{2} < \ldots{})$ be an infinite sequence of positive integers. For $h \ge 2$ integer let $R_{h}(\mathcal{A}, n)$ denote the number of solutions of the equation

\begin{equation}\label{1}
a_{i_{1}} + a_{i_{2}} + \dots + a_{i_{h}} = n, \hspace*{3mm} a_{i_{1}} \in \mathcal{A}, \dots, a_{i_{h}} \in \mathcal{A}, \hspace*{3mm} a_{i_{1}} \le a_{i_{2}} \le \dots{} \le a_{i_{h}},
\end{equation}
\noindent where $n \in \mathbb{N}$. A (finite or infinite) set $\mathcal{A}$ of positive integers is said to be a Sidon set if all the sums $a + b$ with $a,b \in \mathcal{A}$, $a \le b$  are distinct. In other words $\mathcal{A}$ is a Sidon
set if for every $n$ positive integer $R_{2}(\mathcal{A}, n) \le 1$. We say a set $\mathcal{A} \subset \mathbb{N}$ is an asymptotic basis of order $h$, if every large enough positive integer $n$ can be represented as the sum of $h$
terms from $\mathcal{A}$, i.e., if there exists a positive integer $n_{0}$ such that $R_{h}(\mathcal{A}, n) > 0$ for $n > n_{0}$. In [3] and [4] P. Erd\H os, A. S\'ark\"ozy and V. T. S\'os asked if there exists a Sidon set which is
an asymptotic basis of order 3. The problem was also appears in [10] (with a typo in it: order 2 is written instead of order 3). It is easy to see [5] that a Sidon set cannot be an asymptotic basis of order 2. A few years ago J. M. Deshouillers and A. Plagne in [2] constructed a Sidon set which is an asymptotic basis of order at most 7. In [8] S. Kiss proved the existence of a
Sidon set which is an asymptotic basis of order 5. In this paper we will improve this result by proving that there exists an
asymptotic basis of order 4 which is a Sidon set by using probabilistic methods.
\begin{thm1}
There exists an asymptotic basis of order 4 which is a Sidon set.
\end{thm1}

\noindent Note that at the same time Javier Cilleruelo [1] has proved a slightly stronger result namely the existence of a Sidon set which is an asymptotic basis of order $3 + \varepsilon$. He obtained his result independently from our work by using other probabilistic methods.
Before we prove the above theorem, we give a short survey of the probabilistic method we are working with.

\section{Probabilistic tools}
The proof of Theorem 1 is based on the probabilistic method due to Erdős and Rényi. There is an excellent summary of this method in the Halberstam - Roth book [6]. We use the notation and terminology of this book. First we give a survey of the probabilistic tools and notations which we use in the proof of Theorem 1. Let $\Omega$ denote the set of strictly increasing sequences of positive integers. In this paper we denote the probability of an event $E$ by $\mathbb{P}(E)$ and the expectation of a random variable $\xi $ by $\mathbb{E} (\xi )$.

\begin{lem1}
Let $\theta _{1}, \theta _{2}, \theta _{3}, \ldots $ be real numbers satisfying
\[
0 \le \theta _{n} \le 1 \hspace*{4mm} (n = 1, 2, \ldots ).
\]
\noindent Then there exists a probability space ($\Omega$, $X$, $\mathbb{P}$) with the following two properties:
\begin{itemize}
\item[(i)] For every $n \in \mathbb{N}$, the event $\mathcal{E}^{(n)} = \{\mathcal{A}: \mathcal{A} \in \Omega, n \in \mathcal{A}\}$ is measurable, and $\mathbb{P}(\mathcal{E}^{(n)}) = \theta _{n}$.
\item[(ii)] The events $\mathcal{E}^{(1)}$, $\mathcal{E}^{(2)}$, ... are independent.
\end{itemize}
\end{lem1}
See Theorem 13. in [6], p. 142. We denote the characteristic function of the event $\mathcal{E}^{(n)}$ by $t_{(\mathcal{A}, n)}$ or we can say the the boolean variable $t_{(\mathcal{A}, n)}$ means that:
\[
t_{(\mathcal{A}, n)} =t_n=
\left\{
\begin{aligned}
1 \textnormal{, if } n \in \mathcal{A} \\
0 \textnormal{, if } n \notin \mathcal{A}.
\end{aligned} \hspace*{3mm}
\right.
\]

\noindent Furthermore, for some $\mathcal{A} = \{a_1, a_2, \ldots{}\} \in \Omega$ we denote the number of solutions of
$a_{i_{1}} + a_{i_{2}} + \ldots + a_{i_{h}} = n$ with $a_{i_{1}}, \ldots , a_{i_{h}} \in \mathcal{A}$, $1 \le a_{i_{1}} < a_{i_{2}} \ldots < a_{i_{h}} < n$ by $r_{h}(\mathcal{A},n)$.
Let
\begin{equation}\label{2}
r_{h}(\mathcal{A},n) = \sum_{\substack{(a_{i_{1}}, a_{i_{2}}, \ldots, a_{i_{h}}) \\ 1 \le a_{i_{1}} < \ldots < a_{i_{h}} < n \\ a_{i_{1}} + a_{i_{2}} + \ldots + a_{i_{h}} = n}}t_{(\mathcal{A}, a_{1})} t_{(\mathcal{A}, a_{2})} \ldots t_{(\mathcal{A}, a_{h})}.
\end{equation}
\noindent Let $r_{h}^{*}(\mathcal{A}, n)$ denote the number of those representations of $n$ in the form (\ref{1}) in which there are at least two  equal terms. Thus we have
\begin{equation}\label{3}
R_{h}(\mathcal{A}, n) = r_{h}(\mathcal{A}, n) + r_{h}^{*}(\mathcal{A},n).
\end{equation}
\noindent It is easy to see from (\ref{2}) that $r_{h}(\mathcal{A}, n)$ is the sum of random variables. However, for $h > 2$ these variables are not independent because the same $t_{(\mathcal{A}, a_{i})}$ may appear in many terms. To overcome this problem we need deeper probabilistic tools. Our proof is based on a method of J. H. Kim and V. H. Vu. We give a short survey of this method. Interested reader can find more details in [7], [11], [12], [13]. Assume that $t_{1}, t_{2}, \ldots , t_{n}$ are independent binary (i.e., all $t_{i}$'s are in $\{0,1\}$) random variables. Consider a polynomial $Y=Y(t_1, \ldots t_n)$ in $t_{1}, t_{2}, \ldots, t_{n}$ with degree $k$ (where the degree of this polynomial equals to the maximum of the sum of the exponents of the monomials). We say a polynomial $Y$ is totally positive if it can be written in the form $Y = \sum_{i}e_{i}\Gamma_{i}$, where the $e_{i}$'s are positive and $\Gamma_{i}$ is a product of some $t_{j}$'s. Furthermore, $Y$ is regular if all of its coefficients are between zero and one. We also say $Y$ is simplified, if all of its monomials are square-free (i.e. do not contain any factor of $t _i ^2$), and homogeneous if all the monomials have the same degree. Thus for instance a boolean polynomial is automatically regular and simplified, though not necessarily homogeneous. Given any multi-index $\underline{\alpha}= \left( \alpha _1, \dots , \alpha _n \right) \in \mathbb{N}^{n}$ , we define the partial derivative $\partial ^{\underline{\alpha}}(Y)$ of $Y$ as
\[ \partial ^{\underline{\alpha} }(Y) = \left( \frac{\partial}{\partial t_1}\right)^{\alpha _1} \cdots \left( \frac{\partial}{\partial t_n}\right)^{\alpha _n} Y(t_1, \ldots t_n),\]
and denote the order of $\underline{\alpha}$ as $|\underline{\alpha} |= \alpha _1+ \dots + \alpha _n$. For any order $d \geq 0$, we denote $\mathbb{E}_{d }(Y) = \text{max}_{\alpha : |\alpha |=d} \mathbb{E} \left( \partial ^{\underline{\alpha} } Y\right)$. Thus for instance $\mathbb{E}_0 (Y)= \mathbb{E}(Y)$ and $\mathbb{E}_d (Y)=0$ if $d$ exceeds the degree of $Y$. We also define $\mathbb{E}_{\geq d}(Y)=\text{max}_{d'\geq d} \mathbb{E} _{d'}(Y)$. The following result is due to Kim and Vu.

\begin{lem2} (J. H. Kim and V. H. Vu)
Let $k \geq 1$ and $Y=Y(t_1, \ldots , t_n)$ be a totally positive polynomial of $n$ independent boolean variables $t_1, \ldots , t_n$. Then there exists a constant $C_k>0$ depending only on $k$ (which is the degree of the polynomial) such that
\[
\mathbb{P}\left( |Y-\mathbb{E}(Y)| \ge C_{k}\lambda ^{k-\frac{1}{2}}\sqrt{\mathbb{E}_{\geq 0}(Y)\mathbb{E}_{\geq 1}(Y)}\right) = O_k
\left( e^{-\frac{\lambda}{4}+(k-1)\log n} \right)
\]
for all $\lambda >0$.
\end{lem2}

\noindent See [7] for the proof. Informally this theorem asserts that when the derivatives of $Y$ are smaller on average than $Y$ itself, and the degree of $Y$ is small, then $Y$ is concentrated around its mean. Finally we need the Borel - Cantelli lemma:
\begin{lem3}(Borel-Cantelli)
Let $X_{1}, X_{2}, \ldots{}$ be a sequence of events in a probability space. If
\[ \sum_{j=1}^{+\infty}\mathbb{P}(X_{j}) < \infty, \]
then with probability 1, at most a finite number of the events $X_{j}$ can occur.
\end{lem3}
\noindent See in [6], p. 135.

\section{Proof of Theorem 1}
Define the sequence $\theta _n$ in Lemma 1 by
\begin{equation}\label{egy}
\theta _n = n^{-\frac{5}{7}},
\end{equation}
that is $\mathbb{P}\left( \{ \mathcal{A}: \mathcal{A} \in \Omega , n \in \mathcal{A}\}\right)=n^{-\frac{5}{7}}$, for $n \in \mathbb{N}$. For a given set $\mathcal{A} \in \Omega$ let the set $\mathcal{B}$ be the following
\begin{equation}\label{ketto}
\mathcal{B}= \left\{ b: b \in \mathcal{A}, \exists a', a'', a''' \in \mathcal{A}: b+a'=a''+a''', a', a'', a''' < b\right\}.
\end{equation}

\noindent Thus $\mathcal{A} \setminus \mathcal{B}$ is a Sidon set. We will prove that $\mathcal{A} \setminus \mathcal{B}$ is an asymptotic basis of order $4$ with probability $1$. This means that there exists integer $N_0$ such that with probability $1$, $r_4(\mathcal{A}\setminus \mathcal{B},n) >0$ for $n \geq N_0$. Since

\[ r_4(\mathcal{A}\setminus \mathcal{B},n)= r_4(\mathcal{A},n)- \left( r_4(\mathcal{A},n)-r_4(\mathcal{A}\setminus \mathcal{B},n)\right) ,\]

\noindent if we get a lower bound for $r_4(\mathcal{A},n)$ and an upper bound for $\left( r_4(\mathcal{A},n)-r_4(\mathcal{A}\setminus \mathcal{B},n)\right)$ then we will have a lower bound for $r_4(\mathcal{A}\setminus \mathcal{B},n)$. So formally we will show that there are positive constants $C_1$ and $N_1$ such that with probability $1$,
\begin{equation}\label{harom}
r_4(\mathcal{A},n) >C_1 n^{\frac{1}{7}},\quad n\geq N_1,
\end{equation}
and there are positive constants $C_2$ and $N_2$ such that with probability $1$,
\begin{equation}\label{negy}
r_4(\mathcal{A},n) - r_4(\mathcal{A}\setminus \mathcal{B},n) < C_2 \left( \log n \right)^{6,5},\quad n\geq N_2.
\end{equation}

\noindent In order to prove (\ref{harom}) and (\ref{negy}) we use Lemma 2.

\noindent We need the following Lemma (see in \cite{nat}, p. 134., Lemma 5.3). For the sake of completeness we sketch the proof.
\begin{lem4}
Let $N \geq 3$, $\alpha , \beta > -1$. Then
\[ \sum_{n=1}^{N-1} n^{\alpha}(N-n)^{\beta} = \Theta _{\alpha ,\beta } \left( N^{\alpha +\beta +1 }\right).\]
\end{lem4}

\noindent \textbf{Proof.}
\[ \sum_{n=1}^{N-1} n^{\alpha}(N-n)^{\beta}= \sum_{1\leq n\leq \frac{N}{2}} n^{\alpha}(N-n)^{\beta}+\sum_{\frac{N}{2} < n < N} n^{\alpha}(N-n)^{\beta} =\]
\[ = \Theta _{\alpha ,\beta } \left( N^{\beta } \sum_{1\leq n\leq \frac{N}{2}} n^{\alpha} \right) + \Theta _{\alpha ,\beta } \left( N^{\alpha} \sum_{\frac{N}{2} < n < N}(N-n)^{\beta} \right) =\]
\[ = \Theta _{\alpha ,\beta } \left( N^{\beta } \int _1^{\frac{N}{2}}x^{\alpha} \mathrm{d}x \right)+ \Theta _{\alpha ,\beta } \left( N^{\alpha} \int _{\frac{N}{2}} ^N x^{\beta } \mathrm{d}x \right) = \Theta _{\alpha ,\beta } \left( N^{\alpha +\beta +1 }\right) . \blacksquare\]

\noindent In the first step we prove (\ref{harom}) by using Lemma 2. To do this, we need the following Lemma.
\begin{lem5} Assume that all of the variables $y_i$'s are different and the $t_{y_i}$'s are random boolean variables.
\begin{enumerate}
\item For every nonzero integer $a_1$ and for every integer $m$
\[\mathbb{E} \left(\sum _{\substack{y_1 \\ a_1y_1=m}} t_{y_1}\right) = O_{a_1} \left( 1\right).\]
\item For every nonzero integers $a_1, a_2$ and for every integer $m$
\[\mathbb{E} \left(\sum _{\substack{(y_1,y_2) \\ a_1y_1+a_2y_2=m}} t_{y_1}t_{y_2}\right) = O_{a_1,a_2} \left( 1\right).\]
\item For every nonzero integers $a_1, a_2, a_3$ and for every integer $m$
\[\mathbb{E} \left(\sum _{\substack{(y_1,y_2,y_3) \\ a_1y_1+a_2y_2+a_3y_3=m}} t_{y_1}t_{y_2}t_{y_3}\right) = O_{a_1,a_2,a_3} \left( 1\right).\]
\end{enumerate}
\end{lem5}

\noindent \textbf{Proof.} (1):
\[ \mathbb{E} \left(\sum _{\substack{y_1 \\ a_1y_1=m}} t_{y_1}\right) =
\begin{cases}
\left( \frac{m}{a_1}\right)^{-\frac{5}{7}} = O_{a_1} \left( 1\right) \quad &\text{if} \quad \frac{m}{a_1} \in \mathbb{Z}^{+} \\
0 \quad &\text{if} \quad \frac{m}{a_1} \not \in \mathbb{Z}^{+}
\end{cases}
\]

\noindent (2): We distinguish two different cases.

\noindent \emph{Case 1.} Assume, that $a_1 >0$, $a_2>0$, thus $m>0$. (Since $y_1$, $y_2$, $a_1$, $a_2$ are nonnegative, therefore $m$ can not be negative, at this case.) Thus applying Lemma 4 we get
\[ \mathbb{E} \left(\sum _{\substack{(y_1,y_2) \\ a_1y_1+a_2y_2=m}} t_{y_1}t_{y_2}\right) = O_{a_1,a_2} \left(\sum _{y_1=1}^{\frac{m}{a_1}} y_1 ^{-\frac{5}{7}} \left( \frac{m-a_1y_1}{a_2}\right)^{-\frac{5}{7}} \right) =\]
\[ O_{a_1,a_2} \left( \sum _{y_1=1}^{\frac{m}{a_1}} (a_1y_1) ^{-\frac{5}{7}} \left( m-a_1y_1\right)^{-\frac{5}{7}}\right) = O_{a_1,a_2} \left( \sum _{y=1}^{m-1} y^{-\frac{5}{7}} (m-y)^{-\frac{5}{7}}\right) =\]
\[ = O_{a_1,a_2} \left( m^{-\frac{3}{7}}\right) =O_{a_1,a_2} \left( 1\right).\]

\noindent \emph{Case 2.} Now assume that $a_1 >0$, $a_2<0$ and $m\geq 0$. (If $m$ is negative, then consider the equation $-a_1y_1-a_2y_2=m$.) We apply Lemma 4 again.
\[ \mathbb{E} \left(\sum _{\substack{(y_1,y_2) \\ a_1y_1+a_2y_2=m}} t_{y_1}t_{y_2}\right) = O_{a_1,a_2} \left( \sum _{y_2=1}^{\infty} y_2 ^{-\frac{5}{7}} \left( \frac{m-a_2y_2}{a_1}\right)^{-\frac{5}{7}} \right) =\]
\[ = O_{a_1,a_2} \left( \sum _{y=1}^{\infty} y ^{-\frac{10}{7}}\right) =O_{a_1,a_2} \left( 1\right).\]
The other cases can be deduced from the aboves. So we leave the details to the reader. ( The case $a_1 < 0$, $a_2 > 0$, either $m\geq 0$ or $m<0$ is almost the same like \emph{Case 2.}, we have to change the role of $a_1$ and $a_2$. The case $a_1 <0$, $a_2<0$, thus $m<0$ is almost the same like \emph{Case 1.}, we have to get $-a_1$ and $-a_2$ instead of $a_1$ and $a_2$. )

\noindent (3): We distinguish three different cases.

\noindent  \emph{Case 1.} Assume, that $a_1 >0$, $a_2>0$, $a_3>0$ thus $m>0$. (Since $y_1$, $y_2$, $y_3$, $a_1$, $a_2$, $a_3$ are nonnegative, therefore $m$ can not be negative, at this case.) Thus applying Lemma 4 we get
\[ \mathbb{E} \left(\sum _{\substack{(y_1,y_2,y_3) \\ a_1y_1+a_2y_2+a_3y_3=m}} t_{y_1}t_{y_2}t_{y_3}\right) = O_{a_1,a_2,a_3} \left( \sum _{y_1=1}^{\frac{m}{a_1}} y_1^{-\frac{5}{7}}\sum _{y_2=1}^{\frac{m-a_1y_1}{a_2}} y_2^{-\frac{5}{7}}\left( \frac{m-a_1y_1-a_2y_2}{a_3}\right)^{-\frac{5}{7}} \right) =\]
\[ =O_{a_1,a_2,a_3} \left( \sum _{y_1=1}^{\frac{m}{a_1}} y_1^{-\frac{5}{7}}\sum _{y_2=1}^{\frac{m-a_1y_1}{a_2}} (a_2y_2)^{-\frac{5}{7}}\left( m-a_1y_1-a_2y_2\right)^{-\frac{5}{7}}\right) =\]
\[ =O_{a_1,a_2,a_3} \left( \sum _{y_1=1}^{\frac{m}{a_1}} y_1^{-\frac{5}{7}}\left( m-a_1y_1\right)^{-\frac{3}{7}} \right) = O_{a_1,a_2,a_3} \left( m^{-\frac{1}{7}}\right) =O_{a_1,a_2,a_3} \left( 1\right) .\]

\noindent  \emph{Case 2.} Now assume that $a_1 >0$, $a_2>0$, $a_3<0$ and $m\geq 0$. Thus applying Lemma 4 again we get
\[ \mathbb{E} \left(\sum _{\substack{(y_1,y_2,y_3) \\ a_1y_1+a_2y_2+a_3y_3=m}} t_{y_1}t_{y_2}t_{y_3}\right) = O_{a_1,a_2,a_3} \left( \sum _{y_3=1}^{\infty} y_3^{-\frac{5}{7}}\sum _{y_1=1}^{\frac{m-a_3y_3}{a_1}} y_1^{-\frac{5}{7}}\left( \frac{m-a_3y_3-a_1y_1}{a_2}\right)^{-\frac{5}{7}} \right) =\]
\[ =O_{a_1,a_2,a_3} \left( \sum _{y_3=1}^{\infty} y_3^{-\frac{5}{7}}\sum _{y_1=1}^{\frac{m-a_3y_3}{a_1}} (a_1y_1)^{-\frac{5}{7}}\left( m-a_3y_3-a_1y_1\right)^{-\frac{5}{7}}\right) = \]
\[ =O_{a_1,a_2,a_3} \left( \sum _{y_3=1}^{\infty} y_3^{-\frac{5}{7}}\left( m-a_3y_3\right)^{-\frac{3}{7}} \right) = O_{a_1,a_2,a_3} \left( \sum _{y=1}^{\infty } y^{-\frac{8}{7}}\right) =O_{a_1,a_2,a_3} \left( 1\right) .\]

\noindent  \emph{Case 3.} Now assume that $a_1 >0$, $a_2>0$, $a_3<0$ and $m < 0$. By Lemma 4 we get
\[ \mathbb{E} \left(\sum _{\substack{(y_1,y_2,y_3) \\ a_1y_1+a_2y_2+a_3y_3=m}} t_{y_1}t_{y_2}t_{y_3}\right) = O_{a_1,a_2,a_3} \left( \sum _{y_3= \lfloor \frac{m}{a_3}\rfloor +1}^{\infty} y_3^{-\frac{5}{7}}\sum _{y_1=1}^{\frac{m-a_3y_3}{a_1}} y_1^{-\frac{5}{7}}\left( \frac{m-a_3y_3-a_1y_1}{a_2}\right)^{-\frac{5}{7}} \right) = \]
\[ =O_{a_1,a_2,a_3} \left( \sum _{y_3=\lfloor \frac{m}{a_3}\rfloor}^{\infty} y_3^{-\frac{5}{7}}\sum _{y_1=1}^{\frac{m-a_3y_3}{a_1}} (a_1y_1)^{-\frac{5}{7}}\left( m-a_3y_3-a_1y_1\right)^{-\frac{5}{7}}\right) = \]
\[ =O_{a_1,a_2,a_3} \left( \sum _{y_3=\lfloor \frac{m}{a_3}\rfloor}^{\infty} y_3^{-\frac{5}{7}}\left( m-a_3y_3\right)^{-\frac{3}{7}} \right) = O_{a_1,a_2,a_3} \left( \sum _{y=1}^{\infty } y^{-\frac{8}{7}}\right) =O_{a_1,a_2,a_3} \left( 1\right) .\]
The other cases can be deduced from the aboves again. We leave the details to the reader. (The case $a_1 >0$, $a_2<0$, $a_3<0$, either $m\geq 0$ or $m<0$ is almost the same as \emph{Case 2.} and as {Case 3}. The case $a_1 <0$, $a_2<0$, $a_3<0$ thus $m<0$ is almost the same as \emph{Case 1.}. In both we have to get $-a_1$, $-a_2$ and $-a_3$ instead of $a_1$, $a_2$ and $a_2$.) $\blacksquare$

\noindent Now we are ready to prove (\ref{harom}). In view of (2) define $Y$ by
\[ Y = r_{4}(\mathcal{A},n) = \sum _{\substack{(x_1, x_2, x_3, x_4)\\ 1\leq x_1 < \dots < x_4\\ x_1+x_2+x_3+x_4=n}} t_{x_1}t_{x_2}t_{x_3}t_{x_4}.\]
We want to use Lemma 2. To do this we have to estimate the expectation of the variable $Y$ and its partial derivatives. Let $\underline{\alpha} = (\alpha _1, \ldots , \alpha _n)$ be a multi-index. In the first step we prove that for $\underline{\alpha }= \underline{0} $
\begin{equation}\label{ot}
\mathbb{E}\left( \partial ^{\underline{\alpha}} Y\right) = \mathbb{E} (Y) = \Theta (n^{\frac{1}{7}}),
\end{equation}
and for $\underline{\alpha } \neq \underline{0} $
\begin{equation}\label{hat}
\mathbb{E}\left( \partial ^{\underline{\alpha}} Y\right) = O(1).
\end{equation}

\noindent Let $\underline{\alpha }= \underline{0} $. By using Lemma 4 we have
\[\mathbb{E}\left( \partial ^{\underline{\alpha}} Y\right) = \mathbb{E} (Y) = \mathbb{E} \left( \sum_{\substack{(x_1,x_2,x_3,x_4)\\ 1 \leq x_1 < \dots < x_4\\x_1+x_2+x_3+x_4=n}} t_{x_1}t_{x_2}t_{x_3}t_{x_4}\right)=\]
\[ = \Theta \left(\sum _{x_1=1}^{n-3} x_1^{-\frac{5}{7}} \sum _{x_2=1}^{n-x_1-2} x_2 ^{-\frac{5}{7}} \sum _{x_3=1}^{n-x_1-x_2-1} x_3^{-\frac{5}{7}} (n-x_1-x_2-x_3)^{-\frac{5}{7}} \right)=\]
\[ = \Theta \left( \sum _{x_1=1}^{n-3} x_1^{-\frac{5}{7}} \sum _{x_2=1}^{n-x_1-2} x_2 ^{-\frac{5}{7}} (n-x_1-x_2)^{-\frac{3}{7}}\right) = \Theta \left( \sum _{x_1=1}^{n-3} x_1^{-\frac{5}{7}} (n-x_1)^{-\frac{1}{7}}\right) = \Theta \left( n^{\frac{1}{7}}\right) , \]
which shows (\ref{ot}).

\noindent Now assume that $\underline{\alpha } = (\alpha _1, \ldots , \alpha _n) \neq \underline{0} $. If there exists an index $i$ such that $\alpha _i \geq 2$ or $|\underline{\alpha}|=\sum _{i=1}^n \alpha _i \geq 5$, then $\partial ^{\underline{\alpha}} Y=0$. It means that in this case $\mathbb{E}(\partial ^{\underline{\alpha}} Y)=0$. So we may assume that for every index $i$, $\alpha _i \leq 1$ and $| \underline{\alpha} | \leq 4$. Let $l_1 < l_2 < \dots < l_t$, for which $\alpha _{l_1}= \ldots = \alpha _{l_t}=1$, $1 \leq t\leq 4$. If $1 \leq \kappa _1 < \ldots < \kappa _{4-t} \leq 4$, $\kappa _j \in \mathbb{N}$, then $\{ 1,2,3,4\} \setminus \{ \kappa _1, \ldots, \kappa _{4-t}\} = \{ r_1, \ldots , r_t\}$, where $r_1< \dots <r_t$ and $x_{r_1}=l_1, \dots , x_{r_t}=l_t$. It means that the variables $x_{r_1}, \dots , x_{r_t} $ occurs in the partial derivative of $Y$ and the variables $x_{\kappa _1}, \dots , x_{\kappa _{4-t}}$ do not. Thus
\[ \mathbb{E}\left( \partial ^{\underline{\alpha}} Y\right) = \sum _{\substack {(\kappa _1, \ldots , \kappa _{4-t})\\ 1 \leq \kappa _1 < \ldots < \kappa _{4-t} \leq 4}} \mathbb{E} \left(  \sum_{\substack{(x_{\kappa _1}, \ldots ,x_{\kappa _{4-t}})\\ \sum_{j=1}^{4-t} x_{\kappa _j}= n- \sum_{j=1}^t l_j}} t_{x_{\kappa _1}}\ldots t_{x_{\kappa _{4-t}}}\right) .\]

\noindent Since $t \geq 1$, thus $4-t \leq 3$. Since the number of the tuples $(\kappa _1, \ldots , \kappa _{4-t})$ is bounded, by Lemma 5 part (2) we get that the expectation is
\[\mathbb{E}\left( \partial ^{\underline{\alpha}} Y\right)=  \sum _{\substack {(\kappa _1, \ldots , \kappa _{4-t})\\ 1 \leq \kappa _1 < \ldots < \kappa _{4-t} \leq 4}} O(1) = O(1),\]
which proves (\ref{hat}).

\noindent From (\ref{ot}) and (\ref{hat}) we get that
\[ \mathbb{E} _{\geq 0}(Y)=\text{max}_{d'\geq 0} \mathbb{E} _{d'}(Y)= \mathbb{E}(Y)=\Theta (n^{\frac{1}{7}}), \]
and
\[ \mathbb{E} _{\geq 1} (Y)= \text{max}_{d'\geq 1} \mathbb{E} _{d'}(Y)=O(1).\]

\noindent Now apply Lemma 2 with $\lambda =20 \log {n}$ and $k=4$. We get that
\begin{equation}\label{het}
\mathbb{P}\left( \left| Y- \mathbb{E}(Y)\right| \geq C_4 (20 \log{n}) ^{3,5} \sqrt{\mathbb{E} _{\geq 0}(Y)\mathbb{E} _{\geq 1}(Y)} \right) = O \left( \frac{1}{n^2}\right) .
\end{equation}

\noindent Thus by (\ref{het}) and Lemma 3 we get that if $n$ is large enough, with probability $1$,
\[ Y= \mathbb{E} (Y) + O \left( (\log{n})^{3,5} \sqrt{\mathbb{E} _{\geq 0}(Y)\mathbb{E} _{\geq 1}(Y)}\right) = \mathbb{E} (Y) + O \left( n^{\frac{1}{14}}(\log{n})^{3,5} \right) =\Theta (n^{\frac{1}{7}}), \]
which means, that (\ref{harom}) holds.

\noindent In the next step we will prove (\ref{negy}), which shows us that the number of those representations in which there is at least one element from the set $\mathcal{B}$ is not too big. Using the definition of the representation functions and the definitions of the sets $\mathcal{A}$ and $\mathcal{B}$ we get
\begin{equation}\label{nyolc}
 r_{4}(\mathcal{A},n)-r_{4}(\mathcal{A}\setminus \mathcal{B},n)= \sum_{\substack{(a_i, a_j, a_k, a_l) \\ a_i < a_j < a_k <  a_l \in \mathcal{A} \\ a_i + a_j + a_k + a_l=n\\ \exists m \in \{ i,j,k,l \},  \exists a_u, a_v, a_z \in \mathcal{A}\\ a_u, a_v, a_z< a_m \\a_m+a_u=a_v+a_z}} t_{a_{i}} t_{ a_{j}} t_{a_{k}} t_{ a_{l}}
 \end{equation}
 To make the analytic calculations easier we estimate (\ref{nyolc}) and we have that
 \begin{multline}\label{nyolc2}
  r_{4}(\mathcal{A},n)-r_{4}(\mathcal{A}\setminus \mathcal{B},n) =
\frac{1}{24} \sum_{\substack{(a_i, a_j, a_k, a_l) \\ a_i, a_j, a_k,  a_l \in \mathcal{A} \quad \text{are distinct} \\ a_i + a_j + a_k + a_l=n\\ \exists m \in \{ i,j,k,l \},  \exists a_u, a_v, a_z \in \mathcal{A}\\ a_u, a_v, a_z< a_m \\a_m+a_u=a_v+a_z}} t_{a_{i}} t_{ a_{j}} t_{a_{k}} t_{ a_{l}} =\\
= \frac{4}{24} \sum_{\substack{(a_i, a_j, a_k, a_l) \\ a_i, a_j, a_k,  a_l \in \mathcal{A} \quad \text{are distinct} \\ a_i + a_j + a_k + a_l=n\\ \exists a_u, a_v, a_z \in \mathcal{A}\\ a_u, a_v, a_z< a_l \\a_l+a_u=a_v+a_z}} t_{a_{i}} t_{ a_{j}} t_{a_{k}} t_{ a_{l}}
\leq \frac{1}{6} \sum_{\substack{(a_i, a_j, a_k, a_l, a_u, a_v, a_z) \\ a_i, a_j, a_k, a_l, a_u, a_v, a_z \in \mathcal{A} \quad \text{are distinct} \\ a_i + a_j + a_k + a_l=n\\ a_u, a_v, a_z< a_l \\a_l+a_u=a_v+a_z}} t_{a_{i}} t_{ a_{j}} t_{a_{k}} t_{ a_{l}}t_{a_{u}} t_{ a_{v}} t_{a_{z}}
\end{multline}
Using the variables $x_i$-s we can write this in the following form
\begin{equation}\label{nyolc3}
 r_{4}(\mathcal{A},n)-r_{4}(\mathcal{A}\setminus \mathcal{B},n) \leq \sum _{\substack{(x_1, x_2, x_3, x_4, x_5, x_6, x_7)\\ x_1+x_2+x_3+x_4=n\\ x_1, x_2, x_3, x_4 \text{are distinct} \\x_4+x_5=x_6+x_7\\ x_5, x_6, x_7 <x_4}} t_{x_1}\ldots t_{x_7}.
\end{equation}

\noindent So this estimation means, that always $x_4$ will be the element in the $4$-tuple $(x_1, x_2, x_3, x_4)$ which hurts the Sidon property. It is easy to see that in the product $t_{x_1}\ldots t_{x_7}$ the $t_{x_i}$ variables are not necessarily independent. So we need to transform (\ref{nyolc3}). For any $7$-tuple $(x_1, x_2, x_3, x_4, x_5, x_6, x_7)$ with condition $x_5, x_6, x_7 <x_4$ let $\{ x_1, x_2, x_3\} \cap \{ x_5, x_6, x_7\} = \{ x_{i_1}, \dots , x_{i_s}\}$, where $1 \leq i_1 < \dots < i_s \leq 3$. Let $\{ x_5, x_6, x_7\} \setminus \{ x_{i_1}, \dots , x_{i_s}\}= \{ x_{h_1}, \dots , x_{h_u}\}$, where $5 \leq h_1 < \dots < h_u \leq 7$ and $u \leq 3-s$. Then for every fix $s$-tuple $(i_1, \dots , i_s)$ there exist $s+u$ tuple $(d_{i_1}, \ldots ,d_{i_s}, b_{h_1}, \ldots , b_{h_u})$ such that we can write the condition $x_4+x_5-x_6-x_7=0$ in the following form:
\begin{equation}\label{nyolc4}
x_4+ d_{i_1}x_{i_1}+ \ldots + d_{i_s}x_{i_s} + b_{h_1}x_{h_1} + \ldots + b_{h_u}x_{h_u}=0,
\end{equation}
where $x_4, x_{i_1}, \ldots , x_{i_s}, x_{h_1}, \ldots ,x_{h_u}$ are different. In (\ref{nyolc4}) $d_{i_j} \neq 0, j= 1, \dots , s$, $b_{h_j} \neq 0, j= 1, \dots, u$, there is only one positive coefficients, which is equal to $1$ and the sum of the negative coefficients is equal to $-2$. Since $t_x^k=t_x$, if $k \geq 1$ then $t_{x_1}\ldots t_{x_7}= t_{x_1}t_{x_2}t_{x_3}t_{x_4}t_{x_{h_1}}\ldots t_{x_{h_u}}$.

\noindent Thus (\ref{nyolc3}) is equal to the following
 \[ \sum _{\substack{(i_1, \ldots , i_s)\\ 1\leq i_1 < \ldots < i_s \leq 3}} \sum _{\substack{(d_{i_1}, \ldots ,d_{i_s}, b_{h_1}, \ldots , b_{h_u})\\ d_{i_j} \neq 0, b_j \neq 0\\ \text{only one term is positive and} =1 \\ \text{the sum of the negative terms is}  -2}} \sum _{\substack{(x_1, x_2, x_3, x_4, x_{h_1}, \ldots , x_{h_u})\\ x_1+x_2+x_3+x_4=n \\ x_4+ d_{i_1}x_{i_1}+ \ldots + d_{i_s}x_{i_s} + b_{h_1}x_{h_1} + \ldots + b_{h_u}x_{h_u}=0\\ x_{i_j}<x_4, j=1, \dots , s; x_{h_{j'}}<x_4, j'= 1, \dots ,u \\x_1, x_2, x_3, x_4, x_{h_1}, \dots , x_{h_u} \text{are distinct}}} t_{x_1}\ldots t_{x_4}t_{x_{h_1}}\ldots t_{x_{h_u}}. \]

\noindent Let the inner sum be
 \[ Y_{d_{i_1}, \ldots ,d_{i_s}, b_{h_1}, \ldots , b_{h_u}} =\sum _{\substack{(x_1, x_2, x_3, x_4, x_{h_1}, \ldots , x_{h_u})\\ x_1+x_2+x_3+x_4=n \\ x_4+ d_{i_1}x_{i_1}+ \ldots + d_{i_s}x_{i_s} + b_{h_1}x_{h_1} + \ldots + b_{h_u}x_{h_u}=0\\ x_{i_j}<x_4, j=1, \dots , s; x_{h_{j'}}<x_4, j'= 1, \dots ,u \\x_1, x_2, x_3, x_4, x_{h_1}, \dots , x_{h_u} \text{are distinct}}} t_{x_1}\ldots t_{x_4}t_{x_{h_1}}\ldots t_{x_{h_u}}. \]

\noindent Since the number of the variables $Y_{d_{i_1}, \ldots ,d_{i_s}, b_{h_1}, \ldots , b_{h_u}}$ is bounded, it is enough to show that for every $s+u$ tuple $(d_{i_1}, \ldots ,d_{i_s}, b_{h_1}, \ldots , b_{h_u})$ with probability $1$,
\begin{equation}\label{kilenc}
Y_{d_{i_1}, \ldots ,d_{i_s}, b_{h_1}, \ldots , b_{h_u}}= O \left( (\log{n})^{6,5}\right) .
\end{equation}
Let's fix an $s+u$ tuple $(d_{i_1}, \ldots ,d_{i_s}, b_{h_1}, \ldots , b_{h_u})$. We will use Lemma 2. So we have to estimate both the expectation of $Y_{d_{i_1}, \ldots ,d_{i_s}, b_{h_1}, \ldots , b_{h_u}}$ and its partial derivatives. First we will show that for every $\underline{\alpha }=(\alpha _1, \ldots , \alpha _n)$,
\begin{equation}\label{tiz}
\mathbb{E}\left( \partial ^{\underline{\alpha }}Y_{d_{i_1}, \ldots ,d_{i_s}, b_{h_1}, \ldots , b_{h_u}} \right) =O(1)
\end{equation}
holds.

\noindent Let $\underline{\alpha }=(\alpha _1, \ldots , \alpha _n)$. If there exists an index $i$ such that $\alpha _i \geq 2$ or $|\underline{ \alpha} |=\sum _{i=1}^n \alpha _i \geq 5+u$, then $\partial ^{\underline{\alpha}} Y_{d_{i_1}, \ldots ,d_{i_s}, b_{h_1}, \ldots , b_{h_u}}=0$. So we may assume that for every index $i$, $\alpha _i \leq 1$ and $ | \underline{\alpha} | \leq 4+u$. Let's fix $\underline{\alpha} =(\alpha _1, \ldots , \alpha _n)$ and let $\underline{\beta } =(\beta _1, \ldots , \beta_n)$, $\underline{\gamma} =(\gamma _1, \ldots , \gamma_n)$, where $\beta_i \in \mathbb{N}$, $\gamma _i \in \mathbb{N}$ and $\underline{\alpha }= \underline{\beta }+ \underline{\gamma }$. Here $\underline{\beta }$ shows the partial derivatives of the variables $x_1, x_2, x_3, x_4$ and $\underline{\gamma }$ shows the derivatives of the variables $x_{h_1}, \dots , x_{h_u}$. So we can write the partial derivatives of $ Y_{d_{i_1}, \ldots ,d_{i_s}, b_{h_1}, \ldots , b_{h_u}}$ in the following form

\begin{multline*}
\partial ^{\underline{\alpha }}Y_{d_{i_1}, \ldots ,d_{i_s}, b_{h_1}, \ldots , b_{h_u}} = \\
=\sum _{\substack{(\underline{\beta}, \underline{\gamma})\\ \underline{\beta}+ \underline{\gamma}=\underline{\alpha}}} \sum _{\substack{(x_1, x_2, x_3, x_4)\\ x_1+ x_2+x_3+x_4=n\\x_1, x_2, x_3, x_4 \text{are distinct}}} \left( \partial ^{\underline{\beta}} t_{x_1} \ldots t_{x_4}\right) \left( \sum _{\substack{(x_{h_1}, \ldots , x_{h_u})\\x_4+ d_{i_1}x_{i_1}+ \ldots + d_{i_s}x_{i_s} + b_{h_1}x_{h_1} + \ldots + b_{h_u}x_{h_u}=0 \\ x_{i_j}<x_4, j=1, \dots , s; x_{h_{j'}}<x_4, j'= 1, \dots ,u \\x_1, x_2, x_3, x_4, x_{h_1}, \dots , x_{h_u} \text{are distinct}}} \partial ^{\underline{\gamma}} t_{x_{h_1}}\ldots t_{x_{h_u}}\right) .
\end{multline*}

\noindent Since the number of pairs $(\underline{\beta}, \underline{\gamma})$ is bounded, it is enough to show that for every fixed pair $(\underline{\beta}, \underline{\gamma})$, where $0 \leq \beta _i \leq 1$, $0 \leq \gamma _i \leq 1$, $ | \underline{\beta} |+|\underline{\gamma} | = \sum _{i=1}^n (\beta _i+\gamma _i) \leq 4+u$,
\begin{equation}\label{tizenketto}
\mathbb{E} \left( \sum _{\substack{(x_1, x_2, x_3, x_4)\\ x_1+ x_2+x_3+x_4=n\\x_1, x_2, x_3, x_4 \text{are distinct}}} \left( \partial ^{\underline{\beta}} t_{x_1} \ldots t_{x_4}\right) \left( \sum _{\substack{(x_{h_1}, \ldots , x_{h_u})\\x_4+ d_{i_1}x_{i_1}+ \ldots + d_{i_s}x_{i_s} + b_{h_1}x_{h_1} + \ldots + b_{h_u}x_{h_u}=0 \\ x_{i_j}<x_4, j=1, \dots , s; x_{h_{j'}}<x_4, j'= 1, \dots ,u\\x_1, x_2, x_3, x_4, x_{h_1}, \dots , x_{h_u} \text{are distinct}}} \partial ^{\underline{\gamma}} t_{x_{h_1}}\ldots t_{x_{h_u}}\right) \right) =O(1)
\end{equation}
holds. Let's fix now a pair $(\underline{\beta}, \underline{\gamma})$. We will show that for every $4$-tuple $(x_1, x_2, x_3, x_4)$
\begin{equation}\label{tizenharom}
\mathbb{E} \left(\sum _{\substack{(x_{h_1}, \ldots , x_{h_u})\\x_4+ d_{i_1}x_{i_1}+ \ldots + d_{i_s}x_{i_s} + b_{h_1}x_{h_1} + \ldots + b_{h_u}x_{h_u}=0 \\ x_{i_j}<x_4, j=1, \dots , s; x_{h_{j'}}<x_4, j'= 1, \dots ,u\\x_1, x_2, x_3, x_4, x_{h_1}, \dots , x_{h_u} \text{are distinct}}} \partial ^{\underline{\gamma}} t_{x_{h_1}}\ldots t_{x_{h_u}}\right) = O(1).
\end{equation}

\noindent Let $| \underline{\gamma} | =t$, $\gamma _{l_1}= \ldots = \gamma _{l_t}=1$, $1 \leq l_1 < \ldots < l_t \leq n$. Let $\{ g_1, \dots , g_t\} \subseteq \{ h_1, \dots h_u\}$ and $\{ j_1, \dots , j_{u-t}\}= \{ h_1, \dots , h_u\} \setminus \{ g_1, \dots , g_t\}$. These sets shows us that the variables $x_{g_1}, \dots , x_{g_t} $ occurs in the partial derivative of $Y_{d_{i_1}, \ldots ,d_{i_s}, b_{h_1}, \ldots , b_{h_u}}$ and the variables $x_{j_1}, \dots ,x_{j_{u-t}}$ do not. Thus we have
\begin{multline*}
\mathbb{E} \left(\sum _{\substack{(x_{h_1}, \ldots , x_{h_u})\\x_4+ d_{i_1}x_{i_1}+ \ldots + d_{i_s}x_{i_s} + b_{h_1}x_{h_1} + \ldots + b_{h_u}x_{h_u}=0 \\ x_{i_j}<x_4, j=1, \dots , s; x_{h_{j'}}<x_4, j'= 1, \dots ,u\\x_1, x_2, x_3, x_4, x_{h_1}, \dots , x_{h_u} \text{are distinct}}} \partial ^{\underline{\gamma}} t_{x_{h_1}}\ldots t_{x_{h_u}}\right) = \\
= \sum_{\substack{ (j_1, \ldots , j_{u-t})\\ \{ j_1, \dots , j_{u-t}\} \subseteq \{ h_1, \dots , h_u\}}} \sum _{\substack{\pi \\ (x_{g_1}, \ldots , x_{g_t})= \pi (l_1, \ldots , l_t)}} \mathbb{E}\left( \sum _{\substack{(x_{j_1}, \ldots , x_{j_{u-t}})\\ \sum _{q=1}^{u-t} b_{j_q}x_{j_q}= -x_4-d_{i_1}x_{i_1}- \ldots -d_{i_s}x_{i_s}- \sum _{q=1}^t b_{g_q}x_{g_q}}} t_{x_{j_1}} \ldots t_{x_{j_{u-t}}} \right) ,
\end{multline*}
where $\pi (l_1, \ldots , l_t)$ denotes the permutations of $(l_1, \ldots , l_t)$.

\noindent Since the number of the $u-t$ tuples $(j_1, \dots ,j_{u-t})$ and the permutations are bounded and since $u-t \leq 3$, thus from Lemma 4 we get that this expectation is bounded i.e.,
\begin{equation}\label{tizennegy}
 \mathbb{E}\left( \sum _{\substack{(x_{j_1}, \ldots , x_{j_{u-t}})\\ \sum _{q=1}^{u-t} b_{j_q}x_{j_q}= -x_4-d_{i_1}x_{i_1}- \ldots d_{i_s}x_{i_s}- \sum _{q=1}^t b_{g_q}x_{g_q}}} t_{x_{j_1}} \ldots t_{x_{j_{u-t}}}\right) = O(1).
\end{equation}

\noindent So from (\ref{tizennegy}) we get that the left hand side of the equation (\ref{tizenketto}) is equal to the following
\begin{multline*}
\sum _{\substack{(x_1, x_2, x_3, x_4)\\ x_1+ x_2+x_3+x_4=n\\x_1, x_2, x_3, x_4 \text{are distinct}}} \mathbb{E}\left( \partial ^{\underline{\beta}} t_{x_1} \ldots t_{x_4}\right) \mathbb{E} \left(\sum _{\substack{(x_{h_1}, \ldots , x_{h_u})\\x_4+ d_{i_1}x_{i_1}+ \ldots + d_{i_s}x_{i_s} + b_{h_1}x_{h_1} + \ldots + b_{h_u}x_{h_u}=0 \\ x_{i_j}<x_4, j=1, \dots , s; x_{h_{j'}}<x_4, j'= 1, \dots ,u\\x_1, x_2, x_3, x_4, x_{h_1}, \dots , x_{h_u} \text{are distinct}}} \partial ^{\underline{\gamma}} t_{x_{h_1}}\ldots t_{x_{h_u}}\right) = \\
= O \left( \sum _{\substack{(x_1, x_2, x_3, x_4)\\ x_1+ x_2+x_3+x_4=n \\x_1, x_2, x_3, x_4 \text{are distinct}}} \mathbb{E}\left( \partial ^{\underline{\beta}} t_{x_1} \ldots t_{x_4} \right)\right) .
\end{multline*}
By (\ref{hat}) we get that if $\underline{\beta} \neq 0$ then this term is equal to $O(1)$. So we may assume that $\underline{\beta} =0$. This means, that we have to prove that
\begin{equation}\label{tizenot}
\mathbb{E} \left( \sum _{\substack{(x_1, x_2, x_3, x_4)\\ x_1+ x_2+x_3+x_4=n\\x_1, x_2, x_3, x_4 \text{are distinct}}} t_{x_1} \ldots t_{x_4} \left( \sum _{\substack{(x_{h_1}, \ldots , x_{h_u})\\x_4+ d_{i_1}x_{i_1}+ \ldots + d_{i_s}x_{i_s} + b_{h_1}x_{h_1} + \ldots + b_{h_u}x_{h_u}=0 \\ x_{i_j}<x_4, j=1, \dots , s; x_{h_{j'}}<x_4, j'= 1, \dots ,u\\x_1, x_2, x_3, x_4, x_{h_1}, \dots , x_{h_u} \text{are distinct}}} \partial ^{\underline{\gamma}} t_{x_{h_1}}\ldots t_{x_{h_u}}\right) \right) =O(1).
\end{equation}

 \noindent If $1 \leq i_1 < \ldots < i_s \leq 3$ then let $\{1,2,3\} \setminus \{i_1, \ldots i_s\}=\{ f_1, \ldots , f_{3-s}\}$. Thus (\ref{tizenot}) is equivalent to the following
 \begin{multline*} \mathbb{E} \left( \sum _{\substack{(x_{i_1}, \ldots , x_{i_s}, x_4, x_{h_1}, \ldots , x_{h_u})\\ x_4+ d_{i_1}x_{i_1}+ \ldots + d_{i_s}x_{i_s} + b_{h_1}x_{h_1} + \ldots + b_{h_u}x_{h_u}=0 \\ x_{i_j}<x_4, j=1, \dots , s; x_{h_{j'}}<x_4, j'= 1, \dots ,u\\x_{i_1}, \dots , x_{i_s}, x_4, x_{h_1}, \dots , x_{h_u} \text{are distinct}}} \partial ^{\underline{\gamma}} t_{x_{i_1}} \ldots t_{x_{i_s}}t_{x_4}t_{x_{h_1}} \ldots t_{x_{h_u}}  \right. \times \\
\times \left. \left( \sum _{\substack{(x_{f_1}, \ldots , x_{f_{3-s}})\\ x_{f_1}+ \ldots +x_{f_{3-s}} = n-x_{i_1}-\ldots - x_{i_s}-x_4 \\x_{i_1}, \dots , x_{i_s}, x_4, x_{f_1}, \dots , x_{f_{3-s}} \text{are distinct}}} t_{x_{f_1}} \ldots t_{x_{f_{3-s}}}\right) \right) =
\end{multline*}
 \begin{multline} \label{tizenhat}
 =\sum _{\substack{(x_{i_1}, \ldots x_{i_s}, x_4, x_{h_1}, \ldots , x_{h_u})\\  x_4+ d_{i_1}x_{i_1}+ \ldots + d_{i_s}x_{i_s} + b_{h_1}x_{h_1} + \ldots + b_{h_u}x_{h_u}=0 \\ x_{i_j}<x_4, j=1, \dots , s; x_{h_{j'}}<x_4, j'= 1, \dots ,u\\x_{i_1}, \dots , x_{i_s}, x_4, x_{h_1}, \dots , x_{h_u} \text{are distinct}}} \left( \mathbb{E}\left( \partial ^{\underline{\gamma}} t_{x_{i_1}} \ldots t_{x_{i_s}}t_{x_4}t_{x_{h_1}} \ldots t_{x_{h_u}} \right) \right. \times \\
\times \left. \mathbb{E}\left( \sum _{\substack{(x_{f_1}, \ldots , x_{f_{3-s}})\\ x_{f_1}+ \ldots +x_{f_{3-s}} = n-x_{i_1}-\ldots - x_{i_s}-x_4\\x_{i_1}, \dots , x_{i_s}, x_4, x_{f_1}, \dots , x_{f_{3-s}} \text{are distinct}}} t_{x_{f_1}} \ldots t_{x_{f_{3-s}}}\right) \right) .
\end{multline}
 By using Lemma 5 we have
 \[ \mathbb{E}\left( \sum _{\substack{(x_{f_1}, \ldots , x_{f_{3-s}})\\ x_{f_1}+ \ldots +x_{f_{3-s}} = n-x_{i_1}-\ldots - x_{i_s}-x_4 \\x_{i_1}, \dots , x_{i_s}, x_4, x_{f_1}, \dots , x_{f_{3-s}} \text{are distinct}}} t_{x_{f_1}} \ldots t_{x_{f_{3-s}}}\right) =O(1).\]
It follows that (\ref{tizenhat}) is equal to
 \[ O \left( \mathbb{E} \left( \partial ^{\underline{\gamma}} \sum _{\substack{(x_{i_1}, \ldots x_{i_s}, x_4, x_{h_1}, \ldots x_{h_u})\\ x_4+ d_{i_1}x_{i_1}+ \ldots + d_{i_s}x_{i_s} + b_{h_1}x_{h_1} + \ldots + b_{h_u}x_{h_u}=0 \\ x_{i_j}<x_4, j=1, \dots , s; x_{h_{j'}}<x_4, j'= 1, \dots ,u \\x_{i_1}, \dots , x_{i_s}, x_4, x_{h_1}, \dots , x_{h_u} \text{are distinct}}} t_{x_{i_1}} \ldots t_{x_{i_s}}t_{x_4}t_{x_{h_1}} \ldots t_{x_{h_u}} \right) \right) .\]

\noindent Let $\underline{\gamma}= (\gamma _1, \ldots , \gamma _n)$, $|\underline{\gamma}|=t$, $\gamma _{l_1}= \ldots = \gamma _{l_t}=1$, $1 \leq l_1 < \ldots < \l_t \leq n$. If $\{ v_1, \ldots , v_t\} \subset \{ h_1, \ldots , h_u\}$, then let $\{ e_1, \ldots , e_{u-t}\}= \{ h_1, \ldots, h_u\} \setminus \{ v_1, \ldots , v_t\}$. Thus we have
\[ \mathbb{E} \left( \partial ^{\underline{\gamma}} \sum _{\substack{(x_{i_1}, \ldots x_{i_s}, x_4, x_{h_1}, \ldots x_{h_u})\\ x_4+ d_{i_1}x_{i_1}+ \ldots + d_{i_s}x_{i_s} + b_{h_1}x_{h_1} + \ldots + b_{h_u}x_{h_u}=0 \\ x_{i_j}<x_4, j=1, \dots , s; x_{h_{j'}}<x_4, j'= 1, \dots ,u}} t_{x_{i_1}} \ldots t_{x_{i_s}}t_{x_4}t_{x_{h_1}} \ldots t_{x_{h_u}} \right) = \]
\begin{equation} \label{tizenhet}
= \sum _{\{ v_1, \ldots v_t\} \subset \{ h_1, \ldots , h_u\}} \sum _{\substack{\pi \\ (x_{v_1}, \ldots ,x_{v_t}) = \pi (l_1, \ldots , l_t)} }\mathbb{E} \left( \sum_{\substack{ (x_{i_1}, \ldots , x_{i_s}, x_4, x_{e_1} \dots x_{e_{u-t}})\\ x_4+ d_{i_1}x_{i_1}+ \ldots + d_{i_s}x_{i_s} + b_{h_1}x_{h_1} + \ldots + b_{h_u}x_{h_u}=0 \\ x_{i_j}<x_4, j=1, \dots , s; x_{h_{j'}}<x_4, j'= 1, \dots ,u\\x_{i_1}, \dots , x_{i_s}, x_4, x_{h_1}, \dots , x_{h_u} \text{are distinct}}} t_{x_{i_1}} \ldots t_{x_{i_s}}t_{x_4}t_{x_{e_1}} \ldots t_{x_{e_{u-t}}}\right) ,
\end{equation}
where $\pi (l_1, \ldots , l_t)$ denotes the permutations of $(l_1, \ldots , l_t)$.

\noindent Since $s+u+1$ is the number of the variables in the equation $x_4+ d_{i_1}x_{i_1}+ \ldots + d_{i_s}x_{i_s} + b_{h_1}x_{h_1} + \ldots + b_{h_u}x_{h_u}=0$, it is clear, that $s+u+1$ is equal to $3$ or $4$. So $s+u+1-t \leq 3$, except if $s+u+1=4$ and $t=0$. If $s+u+1-t \leq 3$ then we may use Lemma 4 again, which implies (\ref{tizenot}). So we can only show the case $s+u+1=4$ and $t=0$, which is equivalent to $\underline{\gamma} =0$. To do this, it remains to prove the following Lemma.

\begin{lem6} The following expectation are bounded.
\begin{enumerate}
\item $\displaystyle{ \mathbb{E} \left( \sum _{\substack{(x_1, x_2, x_3, x_4)\\ x_1+x_2+x_3+x_4=n}} t_{x_1}t_{x_2}t_{x_3}t_{x_4} \left( \sum _{\substack{(x_5, x_6, x_7)\\ x_4+x_5=x_6+x_7}} t_{x_5}t_{x_6}t_{x_7}\right)\right)} =O(1)$
\item $\displaystyle{\mathbb{E} \left( \sum _{\substack{(x_1, x_2, x_3, x_4)\\ x_1+x_2+x_3+x_4=n}} t_{x_1}t_{x_2}t_{x_3}t_{x_4} \left( \sum _{\substack{(x_5, x_6)\\ x_4+x_3=x_5+x_6}} t_{x_5}t_{x_6}\right)\right)}=O(1)$
\item $\displaystyle{\mathbb{E} \left( \sum _{\substack{(x_1, x_2, x_3, x_4)\\ x_1+x_2+x_3+x_4=n}} t_{x_1}t_{x_2}t_{x_3}t_{x_4} \left( \sum _{\substack{(x_5, x_6)\\ x_4+x_5=x_3+x_6}} t_{x_5}t_{x_6}\right)\right)}=O(1)$
\item $\displaystyle{\mathbb{E} \left( \sum _{\substack{(x_1, x_2, x_3, x_4)\\ x_1+x_2+x_3+x_4=n}} t_{x_1}t_{x_2}t_{x_3}t_{x_4} \left( \sum _{\substack{x_5\\ x_4+x_3=x_2+x_5}} t_{x_5}\right)\right)}=O(1)$
\item $\displaystyle{\mathbb{E} \left( \sum _{\substack{(x_1, x_2, x_3, x_4)\\ x_1+x_2+x_3+x_4=n}} t_{x_1}t_{x_2}t_{x_3}t_{x_4} \left( \sum _{\substack{x_5\\ x_4+x_5=x_2+x_3}} t_{x_5}\right)\right)}=O(1)$
\item $\displaystyle{\mathbb{E} \left( \sum _{\substack{(x_1, x_2, x_3, x_4)\\ x_1+x_2+x_3+x_4=n \\ x_4+x_3=x_2+x_1}} t_{x_1}t_{x_2}t_{x_3}t_{x_4}\right)}=O(1)$
\end{enumerate}
\end{lem6}

\textbf{Proof.} (1): Using the definition of the variables $t_i$'s get that the order of this expectation is
\begin{equation}\label{tizennyolc}
O \left( \sum _{x_1=1}^n x_1 ^{-\frac{5}{7}} \sum _{x_2=1}^{n-x_1} x_2 ^{-\frac{5}{7}} \sum _{x_3=1}^{n-x_1-x_2-1} x_3 ^{-\frac{5}{7}} (n-x_1-x_2-x_3) ^{-\frac{5}{7}} \sum _{x_5=1}^{x_4-1} x_5 ^{-\frac{5}{7}}\sum _{x_6=1}^{x_4+x_5-1} x_6 ^{-\frac{5}{7}} (x_4+x_5-x_6) ^{-\frac{5}{7}} \right) .
\end{equation}
Applying Lemma 4 to the last sum  it follows that (\ref{tizennyolc}) is equivalent to the following
\begin{equation} \label{tizenkilenc}
O \left( \sum _{x_1=1}^n x_1 ^{-\frac{5}{7}} \sum _{x_2=1}^{n-x_1} x_2 ^{-\frac{5}{7}} \sum _{x_3=1}^{n-x_1-x_2-1} x_3 ^{-\frac{5}{7}} (n-x_1-x_2-x_3) ^{-\frac{5}{7}} \sum _{x_5=1}^{x_4-1} x_5 ^{-\frac{5}{7}}(x_4+x_5) ^{-\frac{3}{7}} \right) .
\end{equation}
Since
\[  \sum _{x_5=1}^{x_4-1} x_5 ^{-\frac{5}{7}}(x_4+x_5) ^{-\frac{3}{7}} = O\left( x_4 ^{-\frac{3}{7}} \sum _{x_5=1}^{x_4-1} x_5 ^{-\frac{5}{7}}\right) = O \left( x_4 ^{-\frac{1}{7}}\right) = O \left( (n-x_1-x_2-x_3) ^{-\frac{1}{7}}\right),\]
it follows that (\ref{tizenkilenc}) is equal to
\begin{multline}\label{husz}
O \left( \sum _{x_1=1}^n x_1 ^{-\frac{5}{7}} \sum _{x_2=1}^{n-x_1} x_2 ^{-\frac{5}{7}} \sum _{x_3=1}^{n-x_1-x_2-1} x_3 ^{-\frac{5}{7}} (n-x_1-x_2-x_3) ^{-\frac{6}{7}} \right) =\\
= O \left( \sum _{x_1=1}^n x_1 ^{-\frac{5}{7}} \sum _{x_2=1}^{n-x_1} x_2 ^{-\frac{5}{7}} (n-x_1-x_2) ^{-\frac{4}{7}}\right) = O \left( \sum _{x_1=1}^n x_1 ^{-\frac{5}{7}} (n-x_1) ^{-\frac{2}{7}}\right) =O(1).
\end{multline}

\noindent (2): Using again the definition of the variables $t_i$'s this expression is equivalent with the following
\begin{equation}\label{huszonegy}
\sum _{x_1=1}^n x_1 ^{-\frac{5}{7}} \sum _{x_2=1}^{n-x_1} x_2 ^{-\frac{5}{7}} \sum _{x_3=1}^{n-x_1-x_2-1} x_3 ^{-\frac{5}{7}} (n-x_1-x_2-x_3) ^{-\frac{5}{7}} \sum _{x_5=1}^{x_4+x_3-1} x_5 ^{-\frac{5}{7}}(x_4+x_3-x_5) ^{-\frac{5}{7}}.
\end{equation}
Since
\[  \sum _{x_5=1}^{x_4+x_3-1} x_5 ^{-\frac{5}{7}}(x_4+x_3-x_5) ^{-\frac{5}{7}} = O\left( (x_4+x_3) ^{-\frac{3}{7}} \right) =O \left( x_4 ^{-\frac{3}{7}}\right)= O \left( x_4 ^{-\frac{1}{7}}\right) = O \left( (n-x_1-x_2-x_3) ^{-\frac{1}{7}}\right),\]
it follows that (\ref{huszonegy}) is equal to
\[ O \left( \sum _{x_1=1}^n x_1 ^{-\frac{5}{7}} \sum _{x_2=1}^{n-x_1} x_2 ^{-\frac{5}{7}} \sum _{x_3=1}^{n-x_1-x_2-1} x_3 ^{-\frac{5}{7}} (n-x_1-x_2-x_3) ^{-\frac{6}{7}} \right).\]
It follows from (\ref{husz}) that this is equal to $O(1)$.\\

\noindent (3): Since $x_3 <x_4$, thus $x_3< \frac{x_3+x_4}{2}= \frac{n-x_1-x_2}{2}$. This expectation is
\begin{equation}\label{huszonketto}
 O \left( \sum _{x_1=1}^n x_1 ^{-\frac{5}{7}} \sum _{x_2=1}^{n-x_1} x_2 ^{-\frac{5}{7}} \sum _{x_3=1}^{\frac{n-x_1-x_2}{2}} x_3 ^{-\frac{5}{7}} (n-x_1-x_2-x_3) ^{-\frac{5}{7}} \sum _{x_5=1}^{x_4-1} x_5 ^{-\frac{5}{7}}(x_4+x_5-x_3) ^{-\frac{5}{7}} \right) ,
\end{equation}
where
\[  \sum _{x_5=1}^{x_4-1} x_5 ^{-\frac{5}{7}}(x_4+x_5-x_3) ^{-\frac{5}{7}} = \sum _{x_5=1}^{x_4-x_3} x_5 ^{-\frac{5}{7}}(x_4+x_5-x_3) ^{-\frac{5}{7}}+ \sum _{x_5=x_4-x_3+1}^{x_4-1} x_5 ^{-\frac{5}{7}}(x_4+x_5-x_3) ^{-\frac{5}{7}}  =\]
\[ = O\left( (x_4-x_3) ^{-\frac{5}{7}}\sum _{x_5=1}^{x_4-x_3} x_5 ^{-\frac{5}{7}} \right)+ O \left( \int _{x_4-x_3}^{\infty}x^{-\frac{10}{7}} \mathrm{d}x \right) =O \left( (x_4-x_3) ^{-\frac{3}{7}}\right) .\]
(\ref{huszonketto}) is equal to
\[ O \left( \sum _{x_1=1}^n x_1 ^{-\frac{5}{7}} \sum _{x_2=1}^{n-x_1} x_2 ^{-\frac{5}{7}} \sum _{x_3=1}^{\frac{n-x_1-x_2}{2}} (2x_3) ^{-\frac{5}{7}} (n-x_1-x_2-2x_3) ^{-\frac{6}{7}} \right).\]
It follows from (\ref{husz}) that this is equal to $O(1)$.\\

\noindent (4): Note, that here $n=x_1+2x_2+x_5$. This expectation is
\begin{equation}\label{huszonketto}
 O \left( \sum _{x_1=1}^n x_1 ^{-\frac{5}{7}} \sum _{x_2=1}^{\frac{n-x_1}{2}} x_2 ^{-\frac{5}{7}} \sum _{x_3=1}^{\frac{n-x_1-x_2}{2}} x_3 ^{-\frac{5}{7}} (n-x_1-x_2-x_3) ^{-\frac{5}{7}} (x_4+x_3-x_2) ^{-\frac{5}{7}} \right) ,
\end{equation}
where
\[ (x_4+x_3-x_2) ^{-\frac{5}{7}} = (n-x_1-x_2-x_2) ^{-\frac{5}{7}} = (n-x_1-2x_2) ^{-\frac{5}{7}}.\]
(\ref{huszonketto}) is equal to
\[ O \left( \sum _{x_1=1}^n x_1 ^{-\frac{5}{7}} \sum _{x_2=1}^{\frac{n-x_1}{2}} x_2 ^{-\frac{5}{7}} (n-x_1-2x_2) ^{-\frac{5}{7}} \sum _{x_3=1}^{\frac{n-x_1-x_2}{2}} x_3 ^{-\frac{5}{7}} (n-x_1-x_2-x_3) ^{-\frac{5}{7}} \right) =\]
\[ = O \left( \sum _{x_1=1}^n x_1 ^{-\frac{5}{7}} \sum _{x_2=1}^{\frac{n-x_1}{2}} x_2 ^{-\frac{5}{7}} (n-x_1-2x_2) ^{-\frac{5}{7}} (n-x_1-x_2) ^{-\frac{3}{7}} \right) = \]
\[ =O \left( \sum _{x_1=1}^n x_1 ^{-\frac{5}{7}} \sum _{x_2=1}^{\frac{n-x_1}{2}} (2x_2) ^{-\frac{5}{7}} (n-x_1-2x_2) ^{-\frac{5}{7}} \right) =\]
\[ = O\left( \sum _{x_1=1}^n x_1 ^{-\frac{5}{7}} (n-x_1) ^{-\frac{3}{7}}\right) = O\left( n ^{-\frac{1}{7}} \right)= O(1).\]

\noindent (5): It is clear that the $n=x_1+x_2+x_3+x_4$, $x_4+x_5=x_2+x_3$ equation-system is equivalent with $n=x_1+2x_4+x_5$, $x_4+x_5=x_2+x_3$. Thus this expectation is
\begin{equation}\label{huszonharom}
O \left( \sum _{x_1=1}^n x_1 ^{-\frac{5}{7}} \sum _{x_4=1}^{\frac{n-x_1}{2}} x_4 ^{-\frac{5}{7}} (n-x_1-2x_4) ^{-\frac{5}{7}} \sum _{x_3=1}^{x_4+x_5-1} x_3 ^{-\frac{5}{7}} (x_4+x_5-x_3) ^{-\frac{5}{7}} \right),
\end{equation}
\[ O \left( \sum _{x_3=1}^{x_4+x_5-1} x_3 ^{-\frac{5}{7}}(x_4+x_5-x_3) ^{-\frac{5}{7}} \right) =  O\left( (x_4+x_5) ^{-\frac{3}{7}} \right) =O(1) .\]
So (\ref{huszonharom}) is equal to
\[ O \left( \sum _{x_1=1}^n x_1 ^{-\frac{5}{7}} \sum _{x_4=1}^{\frac{n-x_1}{2}} (2x_4) ^{-\frac{5}{7}} (n-x_1-2x_4) ^{-\frac{5}{7}} \right) =\]
\[ = O \left( \sum _{x_1=1}^n x_1 ^{-\frac{5}{7}}(n-x_1) ^{-\frac{3}{7}} \right) = O\left( n ^{-\frac{1}{7}} \right) = O(1).\]

\noindent (6): It is clear that the $n=x_1+x_2+x_3+x_4$, $x_1+x_2=x_3+x_4$ equation-system is equivalent with $\frac{n}{2}=x_1+x_2=x_3+x_4$. Thus this expectation is
\begin{multline}\label{huszonharom}
O \left( \sum _{x_1=1}^{\frac{n}{2}} x_1 ^{-\frac{5}{7}} \left( \frac{n}{2}-x_1\right) ^{-\frac{5}{7}} \sum _{x_3=1}^{\frac{n}{2}} x_3 ^{-\frac{5}{7}} \left( \frac{n}{2}-x_1\right) ^{-\frac{5}{7}} \right) = \\
= O \left( \sum _{x_1=1}^{\frac{n}{2}} x_1 ^{-\frac{5}{7}} \left( \frac{n}{2}-x_1\right) ^{-\frac{5}{7}} n ^{-\frac{3}{7}}\right) = O\left( n ^{-\frac{6}{7}} \right) =O(1). \blacksquare
\end{multline}

\noindent So from (\ref{tiz}) we get that
\[ \mathbb{E}_{\geq 0}(Y_{d_{i_1}, \ldots ,d_{i_s}, b_{h_1}, \ldots , b_{h_u}}) =\text{max}_{d'\geq 0} \mathbb{E} _{d'}(Y_{d_{i_1}, \ldots ,d_{i_s}, b_{h_1}, \ldots , b_{h_u}}) =O(1),\]
and
\[ \mathbb{E}_{\geq 1}(Y_{d_{i_1}, \ldots ,d_{i_s}, b_{h_1}, \ldots , b_{h_u}}) =\text{max}_{d'\geq 1} \mathbb{E} _{d'}(Y_{d_{i_1}, \ldots ,d_{i_s}, b_{h_1}, \ldots , b_{h_u}})=O(1).\]
 \noindent In Lemma 2, we have $k\leq 7$ and $\lambda = 32 \log n$. Thus

\begin{multline}\label{tizenegy}
\mathbb{P}\left( \left| Y_{d_{i_1}, \ldots ,d_{i_s}, b_{h_1}, \ldots , b_{h_u}}- \mathbb{E}(Y_{d_{i_1}, \ldots ,d_{i_s}, b_{h_1}, \ldots , b_{h_u}})\right| \right. \geq \\
\geq \left. C_k (32 \log{n}) ^{k-\frac{1}{2}} \sqrt{\mathbb{E} _{\geq 0}(Y_{d_{i_1}, \ldots ,d_{i_s}, b_{h_1}, \ldots , b_{h_u}})\mathbb{E} _{\geq 1}(Y_{d_{i_1}, \ldots ,d_{i_s}, b_{h_1}, \ldots , b_{h_u}})}\right) = \\
=O_k \left( e^{-8 \log {n}+(k-1)\log {n}}\right) = O_k \left( \frac{1}{n^2}\right) .
\end{multline}
Thus by (\ref{tiz}) and Lemma 2 we get that with probability $1$,
\begin{multline*}
 Y_{d_{i_1}, \ldots ,d_{i_s}, b_{h_1}, \ldots , b_{h_u}}= \mathbb{E} (Y_{d_{i_1}, \ldots ,d_{i_s}, b_{h_1}, \ldots , b_{h_u}}) + \\
+ O \left( (\log{n})^{k-\frac{1}{2}} \sqrt{\mathbb{E} _{\geq 0}(Y_{d_{i_1}, \ldots ,d_{i_s}, b_{h_1}, \ldots , b_{h_u}})\mathbb{E} _{\geq 1}(Y_{d_{i_1}, \ldots ,d_{i_s}, b_{h_1}, \ldots , b_{h_u}})}\right) =\\
= O \left( ( \log{n})^{6,5} \right) ,
\end{multline*}
which shows (\ref{tiz}), which proves (\ref{negy}). And this completes the proof.$\blacksquare$

\end{document}